\documentclass[12pt,a4paper]{amsart} 
\usepackage[utf8]{inputenc} 
\usepackage[top=3cm, bottom=3.5cm, left=2.5cm, right=2.5cm]{geometry}
\usepackage{amsmath,amsthm,amsfonts,amstext,amssymb}
\usepackage{mathrsfs}
\usepackage{enumerate}
\usepackage{hyperref}
\usepackage{dsfont}
\usepackage{color}
\usepackage{graphicx}
\usepackage{xcolor}
\hypersetup{
    colorlinks,
    linkcolor={blue!50!black},
    citecolor={blue!50!black},
    urlcolor={blue!80!black}
}

\def\R{{\mathbb{R}}}

\def\Z{{\mathbb{Z}}}

\newtheorem{theorem}{Theorem}[section]

\theoremstyle{definition}

\newtheorem{remark}[theorem]{Remark}

\numberwithin{equation}{section}

\begin{document}

\title[On the maximum visibility in a ball]{On the maximum visibility in a ball through the vacant set of Poissonian obstacles}
\author{Yingxin Mu}
\address{
  Yingxin Mu,
  University of Leipzig, Institute of Mathematics,
  Augustusplatz 10, 04109 Leipzig, Germany.
}
\email{yingxin.mu@uni-leipzig.de}

\author{Artem Sapozhnikov}
\address{
  Artem Sapozhnikov,
  University of Leipzig, Institute of Mathematics,
  Augustusplatz 10, 04109 Leipzig, Germany.
}
\email{artem.sapozhnikov@math.uni-leipzig.de}

\begin{abstract}
We study the maximum visibility in a ball inside the vacant set of three obstacle models in $\R^d$ with slow decay of spatial correlations and disparate obstacle geometries: Poisson Boolean models with general i.i.d.\ radii distributions, Poisson cylinders and Brownian interlacements. Let $M_r$ be the maximum distance between points $x$ and $y$ in the ball $B(r)$ such that $x$ is visible from $y$. We prove that $M_r$ divided by $q_r$ converges in probability to an explicit model dependent constant, where $q_r=\log r$, except for the Brownian interlacements in dimension $d=3$, where $q_r = \log r\log\log r$. 
\end{abstract}

\maketitle

\section{Introduction}

Let $\mathcal C$ be a random closed subset of $\R^d$, viewed as an opaque obstacle. 
We say that $x$ is \emph{visible} from $y$ (through the vacant set of $\mathcal C$) if the line segment $[x,y]$ does not intersect $\mathcal C$. 
The problem of visibility in random fields of obscuring elements is of interest in various applications, see e.g.\ \cite{Zacks-visibility}. 
The first mathematical study of visibility goes back to Pólya \cite{Polya-visibility}. The probability of visibility to a large distance was studied in \cite{Calka-visibility} for the Poisson-Boolean model with balls of fixed radius and in \cite{ET-visibility} for the Brownian interlacements. A framework to study visibility to a large distance through the vacant set of Poissonian obstacles was developed in \cite{MS-Visibility-AIHP,MS-Visibility-ECP,MS-Visibility-N}, where, in particular, sharp bounds on the probability of (multi-hop) visibility to a large distance were obtained for the Poisson-Boolean models with arbitrary i.i.d.\ random radii, the Poisson cylinders and the Brownian interlacements. Visibility has also been studied in hyperbolic spaces, see e.g.\ \cite{BJST-visibility-H,TC-visibility-H}.

 \smallskip

In this paper, we are interested in the maximal visibility within a large ball through the vacant set of the Poisson-Boolean models with arbitrary i.i.d.\ random radii, the Poisson cylinders and the Brownian interlacements---Poisson clouds of, respectively, Euclidean balls, doubly-infinite cylinders, and doubly-infinite Wiener sausages; we refer to \cite{MR-Book,TW-cylinders,Sznitman-BI} or to \cite[Section~2]{MS-Visibility-N} for their precise definitions.
For $r>0$, let
\[
M_r=\sup\big\{||x-y||\,:\, x,y\in B(r),\text{ $x$ is visible from $y$}\big\}.
\]
The main result of this paper is that, for each of the three models, $M_r/q_r$ converges in probability to an explicit constant, where $q_r=\log r$, except for the Brownian interlacements in dimension $d=3$, where $q_r = \log r\log\log r$. Our approach is general and can be applied to other models with Poissonian obstacles.

\smallskip

For $n\geq 0$, let $\kappa_n  = \frac{\pi^{n/2}}{\Gamma(\frac n2 + 1)}$ (for $n\geq 1$, $\kappa_n$ is the volume of the unit ball in $\R^n$). 
\begin{itemize}
\item
For $d\geq 2$ and a probability measure $\mathsf Q$ on $\R_+$, let  
\[
\varkappa_{\mathrm{BM}} = \kappa_{d-1}\mathsf E[\varrho^{d-1}],
\]
where $\varrho$ is a random variable with distribution $\mathsf Q$. 
\item
For $d\geq 2$ and $\rho>0$, let 
\[
\varkappa_{\mathrm{PC}} =\kappa_{d-2}\rho^{d-2}\frac{\mathrm{Beta}(\frac12,\frac d2)}{\mathrm{Beta}(\frac12,\frac{d-1}2)}.
\]
\item
For $d\geq 3$ and $\rho>0$, 
let 
\[
\varkappa_{\mathrm{BI}} = \rho^{d-3}\left\{\begin{array}{ll} \pi & d=3\\ 2\pi^{\frac{d-1}2}/\Gamma(\frac{d-3}{2}) & d\geq 4\end{array}\right.
\]
\end{itemize}
\begin{theorem}\label{thm:main}
Let $\alpha>0$ and $\rho>0$. Let $\mathsf Q$ be a probability measure on $\R_+$ with 
\begin{equation}\label{eq:BM-expected-volume}
\mathsf E[\varrho^d]<\infty, 
\end{equation}
where $\varrho$ is a random variable with distribution $\mathsf Q$.\footnote{Condition \eqref{eq:BM-expected-volume} is necessary and sufficient for the Poisson-Boolean model not covering the whole $\R^d$ almost surely, see e.g.\ \cite[Proposition~3.1]{MR-Book}.}
Let $\mathcal C$ be either the Poisson-Boolean model with intensity $\alpha$ and radius distribution $\mathsf Q$ ($d\geq 2$) or the Poisson cylinders at level $\alpha$ with radius $\rho$ ($d\geq 2$) or the Brownian interlacements at level $\alpha$ with radius $\rho$ ($d\geq 3$). Let $q_r = \log r\log\log r$ for the Brownian interlacements in dimenison $d=3$ and, otherwise, let $q_r = \log r$. Then 
\[
\frac{M_r}{q_r} \stackrel{P}\longrightarrow \frac{d}{\alpha\varkappa},\quad\text{as $r\to\infty$},
\]
where $\varkappa$ equals to 
$\varkappa_{\mathrm{BM}}$ for the Poisson-Boolean models, 
$\varkappa_{\mathrm{PC}}$ for the Poisson cylinders, and 
$\varkappa_{\mathrm{BI}}$ for the Brownian interlacements. 
\end{theorem}
\begin{remark}
In fact, we prove that 
\begin{itemize}
\item[(a)]
for any $\gamma>\frac{d}{\alpha\varkappa}$, there exists $\epsilon = \epsilon(d, \mathrm{law}(\mathcal C), \gamma)>0$ such that 
\[
\mathsf P[M_r>\gamma q_r] = O(r^{-\epsilon}),
\]
\item[(b)]
for any $\gamma<\frac{d}{\alpha\varkappa}$, there exists $\epsilon = \epsilon(d, \mathrm{law}(\mathcal C), \gamma)>0$ such that 
\[
\mathsf P[M_r<\gamma q_r] = O(r^{-\epsilon}) + O(\mathsf E\big[\varrho^d\min(\varrho/r,1)^d\big]), 
\]
where the second term is non-zero only for the Poisson-Boolean models. Note that if $\mathsf E[\varrho^{d+s}]<\infty$ for some $s\in(0,d]$, then 
\[
\mathsf E\big[\varrho^d\min(\varrho/r,1)^d\big]  \leq
\mathsf E\big[\varrho^d\min(\varrho^s/r^s,1)\big]\leq 
\mathsf E[\varrho^{d+s}]/r^s = O(r^{-s}),
\]
thus the second term decays subpolynomially only if $\mathsf E[\varrho^{d+s}]=\infty$ for all $s>0$. 
\end{itemize}
\end{remark}

\section{Notation and preliminaries}

Let $x\in\R^d$, $\rho>0$ and $K\subset \R^d$. We denote by $B(x,\rho)$ the closed Euclidean ball at $x$ with radius $\rho$ and write $B(\rho)$ for $B(0,\rho)$. We denote the closed $\rho$-neighborhood of $K$ by $B(K,\rho$), that is $B(K,\rho) = \bigcup_{x\in K}B(x,\rho)$. The $k$-dimensional Lebesgue measure is denoted by $\lambda_k$. 
We write $e_i$ for the $i$-th vector in the canonical orthonormal basis of $\R^d$ and $x(i)$ for the $i$-th coordinate of $x$ in the canonical basis. 

\smallskip

For any $d\geq 2$, the volume of the closed $\rho$-neighborhood of the line segment $[0,x]$ equals
\[
\lambda_d\big(B([0,x],\rho)\big) = \kappa_d\rho^d + \kappa_{d-1}\rho^{d-1}\|x\|;
\]
in particular, $\varkappa_{\mathrm{BM}}$, defined in the introduction, equals 
\begin{equation}\label{eq:varkappa-BM}
\varkappa_{\mathrm{BM}} = 
\lim\limits_{\|x\|\to\infty}\frac{1}{\|x\|} 
\mathsf E\big[\lambda_d\big(B([0,x],\varrho)\big)\big].
\end{equation}
Let $\xi$ be the orthogonal projection of a uniformly distributed unit vector in $\R^d$ on the hyperplane $\{z\in\R^d\,:\,z(1)=0\}$, then 
\[
\mathsf E\big[\|\xi\|\big] = \frac{\mathrm{Beta}(\frac12,\frac d2)}{\mathrm{Beta}(\frac12,\frac{d-1}2)},
\]
see e.g.\ \cite[6.2.1]{AS:math-functions} (and the proof of (4.1) in \cite{MS-Visibility-ECP}); in particular, $\varkappa_{\mathrm{PC}}$, defined in the introduction, equals 
\begin{equation}\label{eq:varkappa-PC}
\varkappa_{\mathrm{PC}} = \kappa_{d-2}\rho^{d-2}\mathsf E[\|\xi\|] 
= \lim\limits_{r\to\infty}\frac{1}{r} 
\mathsf E\big[\lambda_{d-1}\big(B([0,r\xi],\rho)\cap\{z\,:\,z(1)=0\}\big)\big].
\end{equation}
For any $d\geq 3$, the capacity of the closed $\rho$-neighborhood of the line segment $[0,x]$ equals
\[
\mathrm{cap}\big(B([0,x],\rho)\big) = \varkappa_d\rho^{d-3}\Big(\tfrac{\|x\|}{\log \|x\|}\mathds{1}_{d=3} + \|x\|\mathds{1}_{d\geq 4}\Big)(1+o(1)),\quad\text{as }x\to\infty,
\]
where $\varkappa_3= \pi$ and $\varkappa_d =2\pi^{\frac{d-1}2}/\Gamma(\frac{d-3}{2})$ for $d\geq 4$, see \cite[Lemma~2.1]{MS-Visibility-ECP}; in particular, $\varkappa_{\mathrm{BI}}$, defined in the introduction, equals 
\begin{equation}\label{eq:varkappa-BI}
\varkappa_{\mathrm{BI}} = 
\lim\limits_{\|x\|\to\infty}
\mathrm{cap}\big(B([0,x],\rho)\big)\Big(\tfrac{\|x\|}{\log \|x\|}\mathds{1}_{d=3} + \|x\|\mathds{1}_{d\geq 4}\Big)^{-1}.
\end{equation}

\smallskip

Let $\mathcal C\subseteq \R^d$ be an obstacle set. For $x_1,x_2\in \R^d$, we say that $x_1$ is visible from $x_2$ (through the vacant set of $\mathcal C$) if the line segment $[x_1,x_2]$ does not intersect $\mathcal C$. 
For sets $A_1,A_2$, we say that $A_1$ is visible from $A_2$ if there exist $x_1\in A_1$ and $x_2\in A_2$ such that $x_1$ is visible from $x_2$. We only consider rotation and translation invariant obstacle sets and define 
\[
f(s) = \mathsf P[\text{$0$ is visible from $se_1$}].
\]
By the definition of the three models, see e.g.\ \cite[(2.3), (2.4), (2.5)]{MS-Visibility-N}, 
\begin{equation}\label{eq:f(s)}
f(s) = e^{-\alpha T(s)},
\end{equation}
where $T(s)$ equals to
\begin{itemize}\itemsep4pt
\item[(a)]
$\mathrm{cap}\big(B([0,se_1],\rho)\big)$ for the Brownian interlacements with radius $\rho$, 
\item[(b)]
$\mathsf E\big[\lambda_{d-1}\big(B([0,s\xi],\rho)\cap\{z\,:\,z(1)=0\}\big)\big]$ for the Poisson cylinders with radius $\rho$, \item[(c)]
$\mathsf E[\lambda_d(B([0,s],\varrho))]$ for the Poisson-Boolean model with radii distribution $\mathsf Q$. 
\end{itemize}

\smallskip

Throughout the paper, we use $c$ and $C$ to denote strictly positive constants that depend only on $d$, $\gamma$ and the distribution of the obstacle set $\mathcal C$. The exact values of these constants are allowed to change at each occurrence, even within the same string of inequalities.

\section{Proof of Theorem~\ref{thm:main}: upper bound}

In this section, we prove that for any $\gamma>\frac{d}{\alpha\varkappa}$, there exists $\epsilon = \epsilon(d, \mathrm{law}(\mathcal C), \gamma)>0$ such that 
\begin{equation}\label{eq:main-upper}
\mathsf P[M_r>\gamma q_r] = O(r^{-\epsilon}),\quad\text{as $r\to\infty$}.
\end{equation}

\smallskip

Recall the definition of the visibility window $\delta_s$ from \cite[Theorem~1.1]{MS-Visibility-AIHP}: (a) $\delta_s = \frac1s$ for the Poisson-Boolean models in dimensions $d\geq 2$, the Poisson cylinders in dimensions $d\geq 3$, and the Brownian interlacements in dimensions $d\geq 4$, (b) $\delta_s  = 1$ for the Poisson cylinders in dimension $d=2$, and (c) $\delta_s = \frac{\log^2 s}{s}$ for the Brownian interlacements in dimension $d=3$. 
In \cite{MS-Visibility-AIHP}, we proved that for all $s>1$, 
\begin{equation}\label{eq:visibility-AIHP}
\mathsf P\big[0\text{ is visible from }B\big(s e_1,\delta_s\big)\big] \leq C 
\mathsf P[0\text{ is visible from }s e_1]. 
\end{equation}
The proof of \eqref{eq:visibility-AIHP} begins by observing that the visibility of $0$ from $B(s e_1, \delta_s)$ through the vacant set of obstacles with radius $\rho$ implies the visibility of $0$ from $se_1$ through the vacant set of obstacles with reduced radius $(\rho- \delta_s)_+$. The same observation holds if the original event is replaced by the event that $B(0,\delta_s)$ is visible from $B(se_1,\delta_s)$. Thus, the proof of \eqref{eq:visibility-AIHP} in \cite{MS-Visibility-AIHP} implies actually the stronger statement that for all $s>1$, 
\begin{equation}\label{eq:visibility-AIHP-stronger}
\mathsf P\big[B(0,\delta_s)\text{ is visible from }B\big(s e_1,\delta_s\big)\big] \leq C 
\mathsf P[0\text{ is visible from }s e_1],
\end{equation}
which we will use in the proof of \eqref{eq:main-upper}.

\smallskip

Let $\gamma>\frac{d}{\alpha\varkappa}$ and define $\beta_r = \frac13\delta_{\gamma q_r}$. Consider a covering of $B(r)$ by $N$ balls of radius $\beta_r$,
\[
\big\{B(x_k,\beta_r),\,\, 1\leq k\leq N\big\},
\]
so that $B(x_k,\frac12\beta_r)\cap B(x_l,\frac12\beta_r) = \emptyset$ for all $k\neq l$; note that $N\leq C(\frac{r}{\beta_r})^d$. 

If $M_r>\gamma q_r$, then there exist balls $B(x_k,\beta_r)$ and $B(x_l,\beta_r)$, such that $\gamma q_r - 2\beta_r\leq\|x_k-x_l\|\leq\gamma q_r + 2\beta_r$ and $B(x_k,\beta_r)$ is visible from $B(x_l,\beta_r)$. Since the number of such pairs of balls is at most $N C(\frac{\gamma q_r + 2\beta_r}{\beta_r})^d\leq C\frac{q_r^dr^d}{\beta_r^{2d}}$, by the rotation and translation invariance, we obtain that 
\begin{eqnarray*}
\mathsf P[M_r>\gamma q_r] 
&\leq &C\frac{q_r^dr^d}{\beta_r^{2d}}\,\mathsf P\big[B(0,\beta_r)\text{ is visible from }B\big((\gamma q_r - 2\beta_r)e_1,\beta_r\big)\big]\\
&\leq &C\frac{q_r^dr^d}{\beta_r^{2d}}\,\mathsf P\big[B(0,\delta_{\gamma q_r})\text{ is visible from }B\big(\gamma q_re_1,\delta_{\gamma q_r}\big)\big]\\
&\stackrel{\eqref{eq:visibility-AIHP-stronger}}\leq &C\frac{q_r^dr^d}{\beta_r^{2d}}\,
\mathsf P[0\text{ is visible from } \gamma q_r e_1]
= C\frac{q_r^dr^d}{\beta_r^{2d}}\,
f(\gamma q_r).
\end{eqnarray*}
By \eqref{eq:f(s)}, the asymptotic formulas for the capacity and the volume \eqref{eq:varkappa-BM}--\eqref{eq:varkappa-BI}, and the definitions of $q_r$ and $\varkappa$, we obtain that 
\[
-\log f(\gamma q_r) = \alpha T(\gamma q_r) = \alpha \varkappa\gamma \log r \big(1+ o(1)\big), \text{ as $r\to\infty$}.
\]
Thus,
\[
\mathsf P[M_r>\gamma q_r] 
\leq r^{d - \alpha\varkappa\gamma + o(1)}, 
\]
which implies \eqref{eq:main-upper}. 
\qed

\section{Proof of Theorem~\ref{thm:main}: lower bound}

In this section, we prove that for any $0<\gamma<\frac{d}{\alpha\varkappa}$, there exists $\epsilon = \epsilon(d, \mathrm{law}(\mathcal C), \gamma)>0$ such that 
\begin{equation}\label{eq:main-lower}
\mathsf P[M_r<\gamma q_r] = O(r^{-\epsilon}) + O(\mathsf E\big[\varrho^d\min(\varrho/r,1)^d\big]), 
\end{equation}
where the second error term is non-zero only for the Poisson-Boolean models. 

\smallskip

Let $\gamma<\frac{d}{\alpha\varkappa}$. Consider 
\[
\{x_k\}_{1\leq k\leq N} = (q_r^3\,\Z^d)\cap B(r-\gamma q_r)
\]
and define line segments $\ell_k = x_k + [0,\gamma q_re_1]$. 
Note that $N\asymp (r/q_r^3)^d$ and $\ell_k\subset B(r)$ for all $k$. If $M_r<\gamma q_r$, then for all $k$, $x_k$ is not visible from $x_k + \gamma q_r e_1$, that is $\ell_k\cap\mathcal C\neq\emptyset$. Let 
\[
X = \sum\limits_{i=1}^N\mathds{1}_{\{\ell_k\cap\mathcal C=\emptyset\}},
\]
then by the Chebyshev inequality, 
\[
\mathsf P[M_r<\gamma q_r] \leq \mathsf P[X=0] \leq 
\mathsf P\Big[|X-\mathsf E[X]| \geq \mathsf E[X]\Big]
\leq \frac{\mathsf E[X^2] - \mathsf E[X]^2}{\mathsf E[X]^2}.
\]
Note that $\mathsf E[X] = N f(\gamma q_r) \asymp (r/q_r^3)^d f(\gamma q_r) = r^{d - \alpha\varkappa\gamma + o(1)}$ by \eqref{eq:varkappa-BM}--\eqref{eq:f(s)} and the definitions of $\varkappa$ and $q_r$; thus, by the choice of $\gamma$, $\mathsf E[X] \geq c r^{\epsilon}$ for some $\epsilon>0$. 
Furthermore, 
\begin{multline*}
\mathsf E[X^2] = \sum\limits_{k,l=1}^N \mathsf P\big[\ell_k\cap\mathcal C=\emptyset, \ell_l\cap\mathcal C=\emptyset\big]\\
\begin{aligned}
&\leq\,\, \mathsf E[X] + \mathsf E[X]^2 + \sum_{1\leq k\neq l\leq N}\Big(\mathsf P\big[\ell_k\cap\mathcal C=\emptyset, \ell_l\cap\mathcal C=\emptyset\big] - \mathsf P[\ell_k\cap\mathcal C=\emptyset]\mathsf P[\ell_l\cap\mathcal C=\emptyset]\Big)\\
&\leq\,\, \mathsf E[X] + \mathsf E[X]^2 + \mathsf E[X]\max_{1\leq k\leq N}\sum_{l\neq k}\Big(\mathsf P\big[\ell_l\cap\mathcal C=\emptyset\,\big|\,\ell_k\cap\mathcal C=\emptyset\big] - \mathsf P[\ell_l\cap\mathcal C=\emptyset]\Big),
\end{aligned}
\end{multline*}
hence 
\begin{equation}\label{eq:main-lower-2}
\mathsf P[M_r<\gamma q_r] \leq Cr^{-\epsilon} + \frac{1}{\mathsf E[X]} \max_{1\leq k\leq N}\sum_{l\neq k}\Big(\mathsf P\big[\ell_l\cap\mathcal C=\emptyset\,\big|\,\ell_k\cap\mathcal C=\emptyset\big] - \mathsf P[\ell_l\cap\mathcal C=\emptyset]\Big).
\end{equation}
Let $\{\mathcal C_s\}_{s\in S}$ be the obstacles in the support of the Poisson point process (balls in the case of the Poisson-Boolean models, doubly-infinite cylinders in the case of the Poisson cylinders, or doubly-infinite Wiener sausages in the case of the Brownian interlacements); in particular, $\mathcal C = \bigcup_{s\in S}\mathcal C_s$. 

Fix $1\leq k\leq N$ and let $\mathcal C'$ be the union of those obstacles $\mathcal C_s$ that intersect $\ell_k$ and $\mathcal C''$ the union of the obstacles $\mathcal C_s$ that do not intersect $\ell_k$. Note that $\mathcal C'$ and $\mathcal C''$ are independent and 
\[
\mathsf P\big[\ell_l\cap\mathcal C=\emptyset\,|\,\ell_k\cap\mathcal C=\emptyset\big]
= \mathsf P\big[\ell_l\cap (\mathcal C'\cup\mathcal C'')=\emptyset\,|\,\mathcal C'=\emptyset\big] = 
\mathsf P[\ell_l\cap \mathcal C''=\emptyset].
\]
Furthermore,
\[
\mathsf P[\ell_l\cap \mathcal C''=\emptyset] 
=
\frac{\mathsf P[\ell_l\cap \mathcal C''=\emptyset, \ell_l\cap\mathcal C'=\emptyset]}{\mathsf P[\ell_l\cap\mathcal C'=\emptyset]}
=
\frac{\mathsf P[\ell_l\cap \mathcal C=\emptyset]}{\mathsf P[\ell_l\cap\mathcal C'=\emptyset]}.
\]
Hence
\[
\mathsf P\big[\ell_l\cap\mathcal C=\emptyset\,|\,\ell_k\cap\mathcal C=\emptyset\big] - \mathsf P[\ell_l\cap\mathcal C=\emptyset]
=\mathsf P[\ell_l\cap\mathcal C=\emptyset]\,\frac{\mathsf P[\ell_l\cap\mathcal C'\neq\emptyset]}{\mathsf P[\ell_l\cap\mathcal C'=\emptyset]}.
\]
Note that 
\begin{eqnarray*}
\mathsf P[\ell_l\cap\mathcal C'\neq\emptyset]
&= &\mathsf P\big[\exists s\in S\,:\,\mathcal C_s\cap \ell_k\neq\emptyset, \mathcal C_s\cap \ell_l\neq\emptyset\big]\\
&\leq &\mathsf P\big[\exists s\in S\,:\,\mathcal C_s\cap B(x_k,\gamma q_r)\neq\emptyset, \mathcal C_s\cap B(x_l,\gamma q_r)\neq\emptyset\big].
\end{eqnarray*}
It is well known that the latter probability is bounded from above by 
\[
C\frac{\mathrm{cap}\big(B(x_k, \gamma q_r + \rho)\big)\mathrm{cap}\big(B(x_l,\gamma q_r + \rho)\big)}{(\|x_k - x_l\| - 2\gamma q_r - 2\rho)^{d-2}} 
\leq Cq_r^{2(d-2)}\|x_k - x_l\|^{2-d}
\]
for the Brownian interlacements (see e.g.\ \cite[(1.45)]{Li-BI}), by 
\[
C(\gamma q_r+\rho)^{2(d-1)}\|x_k - x_l\|^{1-d}\\
\leq Cq_r^{2(d-1)}\|x_k - x_l\|^{1-d}
\]
for the Poisson cylinders (see e.g.\ \cite[Lemmas~3.1 and 3.3]{TW-cylinders}), and by 
\begin{multline*}
\mathsf P\big[\exists s\in S\,:\, \mathcal C_s\cap B(x_k,\gamma q_r)\neq \emptyset, \mathrm{diam}(\mathcal C_s)\geq \|x_k-x_l\|-2\gamma q_r\big]\\
\leq C\mathsf E\big[\varrho^d\,;\,\varrho\geq \frac12\|x_k-x_l\| - \gamma q_r\big]\leq C\mathsf E\big[\varrho^d\,;\,\varrho\geq \frac14\|x_k-x_l\|\big]
\end{multline*}
for the Poisson-Boolean models (see e.g.\ \cite[Lemmas~3.4 and 3.5]{Gouere08}). 

Furthermore, since $\|x_k-x_l\|\geq q_r^3$ for all $l\neq k$, there exists $c<1$ such that 
\[
\mathsf P[\ell_l\cap\mathcal C'\neq\emptyset] \leq c<1,\text{ for all $l\neq k$}.
\]
All in all, for each $k$, 
\[
\sum\limits_{l\neq k}\Big(\mathsf P\big[\ell_l\cap\mathcal C=\emptyset\,|\,\ell_k\cap\mathcal C=\emptyset\big] - \mathsf P[\ell_l\cap\mathcal C=\emptyset]\Big) \leq 
C\,\frac{f(\gamma q_r)}{1-c} \sum\limits_{l\neq k}E_r(x_k,x_l),
\]
where $E_r(x_k,x_l)$ equals to $q_r^{2(d-2)}\|x_k - x_l\|^{2-d}$ for the Brownian interlacements, $q_r^{2(d-1)}\|x_k - x_l\|^{1-d}$ for the Poisson cylinders, and $\mathsf E\big[\varrho^d\,;\,\varrho\geq \frac14\|x_k-x_l\|\big]$ for the Poisson-Boolean models. Note that 
\[
\sum\limits_{l\neq k}\|x_k - x_l\|^{2-d} \leq C \Big(\frac{r}{q_r^3}\Big)^2\quad\text{resp.}\quad
\sum\limits_{l\neq k}\|x_k - x_l\|^{1-d} \leq C \frac{r}{q_r^3}.
\]
Thus, for both the Brownian interlacements in dimensions $d\geq 3$ and the Poisson cylinders in dimensions $d\geq 2$, 
\[
\frac{1}{\mathsf E[X]} \max_{1\leq k\leq N}\sum_{l\neq k}\Big(\mathsf P\big[\ell_l\cap\mathcal C=\emptyset\,\big|\,\ell_k\cap\mathcal C=\emptyset\big] - \mathsf P[\ell_l\cap\mathcal C=\emptyset]\Big) \leq C r^{-1} q_r^{3d},
\]
and \eqref{eq:main-lower} follows from \eqref{eq:main-lower-2}. Now, 
\[
\sum\limits_{l\neq k}\mathsf E\big[\varrho^d\,;\,\varrho\geq \frac14\|x_k-x_l\|\big] \leq C\mathsf E\Big[\varrho^d \Big(\frac{\min(\varrho,r)}{q_r^3}\Big)^d\Big],
\]
hence for the Poisson-Boolean models, 
\[
\frac{1}{\mathsf E[X]} \max_{1\leq k\leq N}\sum_{l\neq k}\Big(\mathsf P\big[\ell_l\cap\mathcal C=\emptyset\,\big|\,\ell_k\cap\mathcal C=\emptyset\big] - \mathsf P[\ell_l\cap\mathcal C=\emptyset]\Big) \leq C\mathsf E\big[\varrho^d\min(\varrho/r,1)^d\big], 
\]
and \eqref{eq:main-lower} follows from \eqref{eq:main-lower-2}. 
\qed

\section*{Acknowledgements}
The research of both authors has been supported by the DFG Priority Program 2265 ``Random Geometric Systems'' (Project number 443849139).

\end{document}